\input amstex
\documentstyle{amsppt}
\NoBlackBoxes
\magnification 1200
\vsize=21truecm
\hsize=16truecm
\define \ov{\overline}

\define \Q{\Bbb Q}
\define \M{\frak M}
\define \Z{\Bbb Z}

\define \sgn{\operatorname{sgn}}
\define \diag{\operatorname{diag}}

\define\BKN{BKN-equation }
\define\CP{CP-system }
\topmatter
\title
          Homological Properties of Graph Manifolds
\endtitle
\author              P. Svetlov
\endauthor
\date draft version, 2001
\enddate
\abstract
We consider the following properties
 of compact oriented irreducible graph-manifolds:
 to contain a $\pi_1$-injective  surface
 (immersed, virtually embedded
 or embedded), be (virtually) fibered over $S^1$,
 and to carry a metric of nonpositive sectional curvature.
 It turns out that all these properties
 can be described from a unified point of view.
 \endabstract
 \endtopmatter

         \head
 \S 0 introduction
 \endhead

 We consider the class ${\frak M}$ of graph-manifolds (cf.\cite{BK1}).
 We say that a 3-manifold $M$ belongs to ${\frak M}$
 if $M$ is a compact closed orientable Haken 3-manifold
 and its canonical decomposition surface (JSJ-surface) $\Cal T$ splits $M$
 into Seifert pieces with orientable base orbifolds
 of negative Euler characteristic. These Seifert pieces
are called {\it maximal blocks}.

 We say that a manifold $M\in\M$ {\it contains} an immersed (resp. embedded)
 surface $S$ if there is
 an immersion (resp. embedding) $g\:S \to M$.
 Such a surface is called $\pi_1$-{\it injective}, if
 the induced homomorphism of fundamental groups
 $g_*:\pi_1(S)\to\pi_1(M)$ is injective.
 Such a surface is called {\it horizontal}
 if $g\:S\to M$ is transversal to
 the Seifert fibers in each maximal block of $M$.
 (Horizontality implies $\pi_1$-injectivity \cite{RW}).

 We say that an immersed $\pi_1$-injective surface $g\:S\to M\in\M$
 is {\it virtually embedded} if
 there exists a finite covering
 $p:\widetilde M\to M$ and a $\pi_1$-injective embedding
 $\tilde{g}\:S\to\widetilde M$ such that $g=p\circ\tilde{g}$.
 Finally we say that a manifold $M\in\M$ is {\it virtually fibered}
 over $S^1$ if $M$ has a finite cover 
 which is fibered over the circle.
 Note, that the fiber of such a bundle has a negative
 Euler characteristic.

 Consider the following properties of $M\in\M$:
  \roster
 \item "(Im)" $M$ contains an immersed
 $\pi_1$-injective surface $S$ with $\chi (S)<0$;

 \item "(HI)" $M$ contains an immersed $\pi_1$-injective horizontal
 surface $S$ with $\chi (S)<0$;
 \item "(F)" $M$ is fibered over $S^1$;

 \item "(E)" $M$ contains  an embedded $\pi_1$-injective surface $S$
 with $\chi (S)<0$;
 \item "(VF)" $M$ is virtually fibered over $S^1$;

 \item "(VE)" $M$ contains a virtually
 embedded $\pi_1$-injective surface $S$  with $\chi (S)<0$;

 \item "(NPC)" $M$ carries a metric of nonpositive sectional curvature
 (NPC-metric),

 \endroster
 where $\chi(S)$ is the Euler characteristic of $S$.

 These properties were intensively studied in the last decade
 by different authors.
 In the paper \cite{LW}, a simple obstruction
 for (virtual) fibration of graph manifolds was described.
 The question about existence of NPC-metrics on graph manifolds
 was treated in the paper \cite{L}.
 In the papers \cite{BK1,2}, an explicit criterion
 which permits to solve whether or not a given
 graph manifold of the class $\M$ carries a NPC-metric
 was obtained.
 In the papers \cite{N1,2}, explicit criteria
 which permit to solve whether or not a given
 graph manifold satisfies one of these properties
 {\bf F}, {\bf VE}, and {\bf Im} were proved;
 an implicit criterion for the property {\bf VF}
 was obtained in \cite{N1} as well.
 In the papers \cite{RW,M}, the question
 whether or not a given $\pi_1$-injective immersion of a surface in
 a graph manifold is a virtual embedment was analyzed.
 If a graph manifold $M$ of the class $\M$ carries
 a NPC-metric then it is virtually fibered over the circle,
 as it was proved in the author work \cite{S}.
 \comment
Note that the last assertion
  is trivial for the graph manifolds
  with non-empty boundary:
  {\it every} such manifold  is
  virtually fibered over the circle  \cite{}
  as well as it admits a NPC-metric \cite{L}.
  On the other hand there is a lot of closed graph manifolds
  which do not admit NPC-metrics or are not virtually fibered.
\endcomment
 So, we can present relationships between these properties
 by the following diagram:
 $$
 \matrix
         E       &\Longrightarrow &   VE   & \Longrightarrow  &    Im       \\
     \Uparrow    &                &\Uparrow&                  &\Updownarrow\\
         F       &\Longrightarrow &   VF   & \Longrightarrow  &    HI      \\
                 &                &\Uparrow&                  &            \\
                 &                &   NPC  &                  &
 \endmatrix  $$
 (the implication {\bf Im}$\Leftrightarrow${\bf HI} follows from \cite{N2}).

 As it follows from S.~Buyalo and V.~Kobel'skii paper \cite{BK1},
 a graph manifold $M\in\M$ admits a NPC-metric iff
 a difference equation over a finite graph
 is solvable. From the other side,
 the same equation was independently
 obtained by W.~Neumann in the paper \cite{N2},
 where $\pi_1$-injective surfaces in graph manifolds
 were investigated. It turns out that the equation
 describes all the features of  graph manifolds
 we have considered above (theorem II below).
 The equation  generalizes the Laplace equation on graphs
 (it contains  long ``covariant'' differences, which
 depends on topological invariants of $M$) \cite{B}.
 We give the equation in 1.4 and call it the
 {\it Buyalo-Ko\-bel'\-skii-Neumann equation} (BKN-equation).

 To learn the properties {\bf Im}--{\bf NPC}
 we also propose an equivalent approach which
 is based on their homological nature
 (theorem I below). This approach calls on
 no topological invariants of manifolds except their
 JSJ-splitting so it is more common than BKN-equation.
 However, \BKN is more suitable for computations.

 \head
 \S 1. Main definitions and results
 \endhead

 \subheading {1.1 Notations}
  Let $M$ be a graph-manifold of the class ${\frak M}$
 and ${\Cal T}$ be its JSJ-surface.
  The surface splits $M$ into a finite set of
 Seifert fibered pieces: $M|{\Cal T}=\{M_v|\,v\in V\}$.
 So $M$ can be obtained by pasting together
 Seifert bundles $M_v,\,v\in V$ along their boundary tori.
 Define the index set $W$ by $\{T_w|\,w\in W\}=\{\partial M_v|\,v\in V\}$
and by definition put
 $\partial v=\{w\in W\,|T_w\subset\partial M_v\}$, so that
 we have $\partial M_v=\{T_w|\,w\in\partial v\}$.
 Denote by $-w$ the element of $W$ such that
 the tori $T_w$ and $T_{-w}$
 are glued together in $M$ providing JSJ-torus $T_{|w|}$.
 These collections $V$ and $W$ form a graph $\Gamma_M(V,W)$, which is
 conjugated to the JSJ-splitting of $M$. Namely,
 $V$ is the vertices set and $W$ is the edges one.
 A vertex $v\in V$ is the initial vertex for an edge $w$ iff $w\in\partial v$.

From now on we fix an orientation of $M$ and
 also we fix some orientations of the
 Seifert fibers in the maximal blocks.
 Let $f_v\in H_1(M_v)$ be the class of the Seifert fiber in $M_v$
 and let
 $f_{w}\in H_1(T_{|w|})$, $w\in\partial v$ be the homological class
 which is defined by $(\iota_w)_*f_w=f_v$,
 where $(\iota_w)_*:H_1(T_{|w|})\to H_1(M_v)$ is the map
 induced by inclusion.

 Let $l\in H^1(M_v)$ be a cohomological class
 in a block $M_v$.
 If $T_w$ is a boundary torus of $M_v$
 and $T_{-w}\subset\partial M_{v'}$
 then we define
 $$
 \langle l,f_{v'}\rangle_{|w|}=
 \iota_w^*l(f_{-w}).
 $$

 The aim of this paper is to present all properties we have introduced
 as modifications of the following one:

 \subheading{1.2 Definition} Let $M$ be a  graph-manifold
 of the class $\M$ and $\Gamma_M(V,W)$ be its graph.
 We say that a nontrivial set $\{l_v \in H^1(M_v)|\,v\in V\}$
of (rational) cohomological classes
   {\it satisfies the compatibility property} (or {\it is a CP-system}
 on $M$), if
\roster
 \item"CP1" $|\langle l_{v'},f_{v}\rangle_{|w|}|\le l_v(f_v)$
for each triple $v,w,v'$ such that $w\in\partial{v}$
 and ${-w}\in\partial{v'}$;

 \item"CP2"
if  $|\langle l_{v'},f_{v}\rangle_{|w|}|=
 l_{v}(f_{v})$ and $|\langle l_{v},f_{v'}\rangle_{|w|}|=
 l_{v'}(f_{v'})$
 then $\iota^*_{-w} l_{v'}=\pm\iota^*_{w}l_{v}$;
 \item"CP3"
           if $l_v(f_v)=0$ then $l_v=0$.

\endroster

 Not all graph manifolds have CP-systems.

 \proclaim{Theorem {I}}
 A graph manifold $M$ of the class
 $\M$ satisfies the property {\bf Im} (resp. {\bf HI}, {\bf F},
 {\bf E}, {\bf VF}, {\bf VE}, {\bf NPC})
 if and only if there exists a CP-system
 on $M$
 (resp. there exists a CP-system $\{\,l_v\,|\,{v\in V}\}$ on $M$ such that
 \roster
 \item "(HI)" $l_{v}(f_{v})>0$;

 \item "(F)"$\langle l_{v'},f_{v}\rangle_e=
 \epsilon_{v'}\epsilon_v l_v(f_v)\ne 0$, where
 $\epsilon:V\to\{\pm 1\}$ -- is a function;

 \item "(E)"
 $|\langle l_{v'},f_{v}\rangle_e|=l_v(f_v)$
;

 \item "(VF)" $l_v(f_v)> 0$ and
 $\langle l_{v'},f_{v}\rangle_e\cdot l_{v'}(f_{v'})=
 \langle l_{v},f_{v'}\rangle_e\cdot l_v(f_v)$;

 \item "(VE)"
 $\langle l_{v'},f_{v}\rangle_e\cdot l_{v'}(f_{v'})=
 \langle l_{v},f_{v'}\rangle_e\cdot l_v(f_v)$;

 \item "(NPC)" $|\langle l_{v'},f_{v}\rangle_e| < l_v(f_v)$
 and $\langle l_{v'},f_{v}\rangle_e\cdot l_{v'}(f_{v'})=
 \langle l_{v},f_{v'}\rangle_e\cdot l_{v}(f_{v})$,

 \endroster
 for each triple $v,w,v'\,(e=|w|)$ as in CP1).
 \endproclaim

 This theorem follows from theorems 2.1, 2.2, 3.1, 3.2,
 4.2, 4.3, 5.1. \qed

 \subheading {Remark} S.~Buyalo gave similar conditions for
  {\bf Im}, {\bf HI}, {\bf E}, and {\bf NPC} (unpublished).
 He analyzed collections of {\it homological} classes
 on JSJ-tori and used the Waldhausen basises which were
 the basic construction in \cite{BK1,2}.

 \subheading {1.3 Invariants of graph manifolds}
 Let $M$ be a manifold of the class $\M$.
 The choosing orientation of $M$ induces orientations on
 maximal blocks $\{M_v|\,v\in V\}$ and hence
 on each torus $T_w,\,w\in W$.
 So, we have an ``area" isomorphism
 $\operatorname{is}_w:\Lambda^2H_1(T_{|w|})\to\Q$.
 By definition, put $a\wedge_w b=\operatorname{is}_w(a\wedge b)$.
 Note that $a\wedge_{-w} b=-a\wedge_{w} b$
 inasmuch choosing orientations
 on $T_w$ and $T_{-w}$ are opposite.

 Now we define the {\it intersection index}
 $b_w\in\Z\setminus\{0\},\,w\in W $ by
 $$
 b_w=f_w\wedge_w f_{-w}.
 $$
 It is clear that $b_w=b_{-w}$.

 Consider  a maximal block $M_v$ in $M\in\M$.
 An element  $z_v\in H_1(\partial M_v)$
 is called {\it adjoint} element for the Seifert fibration
 of $M_v$, if $f_w\wedge_wp_w(z_v)=1$
 for each $w\in\partial v$, where $p_w\:H_1(\partial M_v)\to H_1(T_{|w|})$
is the canonical projection ($\partial M_v=\bigcup_{w\in\partial v}T_{|w|}$).
 For each adjoint element $z_v\in H_1(\partial M_v)$
 we have $\iota_*(z_v)=-e_{M_v}(z_v)\cdot f_v$,
 where $e_{M_v}(z_v)$ is a rational number \cite{N2}.
 The ratio $e_{M_v}(z_v)$  is called {\it the Euler number}
 of the Seifert fibration  of $M_v$ with respect to the
 linear foliation of $\partial M_v$
 induced by $z_v\in H_1(\partial M_v)$.
 This number was introduced in \cite{NR}.

\proclaim{1.3.1 Lemma \cite{N2}}
  Let $M_v$  be a maximal Seifert block
in a graph manifold $M$ and let $p\:S\to M_v$
be a proper $\pi_1$-injective immersion of an oriented surface
with negative Euler characteristic. Consider the (oriented) boundary
of $S$ and let $\{C_{wi},\,i=1,\cdots, k_w\}$ be its components on $T_w$.
Then for each adjoint element $z_v\in H^1(\partial M_v)$
we have
          $$
  \sum_{i=1}^{k_w}f_w\wedge_{w}[C_{wi}]=a\ne 0 \quad
        \text{for each} \;T_w\subset\partial M_v.\tag{1}
        $$

$$
            \sum_{w,i}p_w(z_v)\wedge_{w}[C_{wi}]=-ae_{M_v}(z_v).\tag{2}
$$
Conversely, suppose that for each boundary component $T_w$,
$w\in\partial v$ a family $C_{w 1},\dots,C_{w k_w }$ of immersed curves
transverse to the Seifert fibres of $M_v$ is given satisfying homology
relations \thetag{1}, \thetag{2}.  Then there exist integers $d_0>0, n_0>0$
so that for any positive integer multiple $d$ of $d_0$ and
$n$ of $n_0$ the family of curves $d(C_{w i})^n$,
$w\in\partial v$, $i =1,\dots,k_w $,
obtained as follows, bounds a connected immersed horizontal surface.
We take $d(C_{w i})^n$ as d copies of
the immersed curve obtained by going $n$ times around the curve $C_{wi}$.
\endproclaim

\proclaim{Lemma 1.3.2}    Let $M_v$  be a maximal Seifert block
in a graph manifold $M$  and for each $w\in\partial v$
a {\it disjoint} family of simple oriented curves
$\{C_{wi},\,i=1,\cdots,k_w\}\subset T_{|w|}$ is given.
If the conditions \thetag{1} and \thetag{2} are
satisfied then the set of curves bounds
an {\it embedded} horizontal surface
in $M_v$. If the number $a$ in \thetag{1} and \thetag{2}
is zero then there is a set of incompressible, boundary
incompressible disjoint embedded vertical annuli in $M_v$
whose boundary is homotopic to the family of curves.
\endproclaim
\demo{Proof}
Consider the classes
$c_w=\sum_{i=1}^{k_w}[C_{wi}]\in H_1(T_{|w|})$.
For each adjoint element $z_v\in H_1(\partial M_v)$ we can write
$c_w=ap_w(z_v)-(p_w(z_v)\wedge_wc_w)f_w$.
Let $c\in H_1(\partial M_v)$ be the class such that
$p_w(c)=c_w$. Then for any $l\in H^1(M_v)$
we have $l(\iota_*c)=0$:
$$\multline
l(\iota_*c)=\sum_{w\in\partial v}\iota_w^*l(c_w)=
\sum_{w\in\partial v}\iota_w^*l(ap_w(z_v)-(p_w(z_v)\wedge_wc_w)f_w)=\\
a\sum_{w\in\partial v}\iota_w^*l(p_w(z_v))-
\sum_{w\in\partial v}(p_w(z_v)\wedge_wc_w)\iota_w^*l(f_w)=
a\sum_{w\in\partial v}l\Bigl((\iota_w)_*p_w(z_v)\Bigr)- \\
l(f_v)\sum_{w\in\partial v}(p_w(z_v)\wedge_wc_w)=
-ae_{M_v}(z_v)l(f_v)+ae_{M_v}(z_v)l(f_v)=0.
\endmultline
$$
It means that $\iota_*c=0$, so there exists an integer class
$s\in H_2(M_v;\partial M_v)$ such that $s\mapsto c$
under the canonical homomorphism
$H_2(M_v;\partial M_v)\to H_1(\partial M_v)$.
The class $s$ can be realized as an embedded surface $S\subset M_v$.

If $a=0$ then all curves are homotopic to the Seifert fibers (by \thetag{1})
and there are even number of curves (by \thetag{2}).
It is obvious that there are required annuli.
\qed\enddemo

 Now we define the {\it charge} $k_v$ of a maximal block
 $M_v$ in a graph manifold $M\in \M$ by
 $k_v=-e_{M_v}(\phi_v)$, where
 $\phi_v$ is the element of $H_1(\partial M_v)$
 such that $p_w(\phi_v)=f_{-w}/b_w$, $w\in\partial v$.

  The charge $k_v$ of a block $M_v\subset M$
  is independent of the orientations of the Seifert fibers
  in blocks of M. It depends only on
  orientation of $M$.
  See, also \cite{BK1}.

The graph $\Gamma_M(V,W)$ providing with the numbers
$\{\,b_w,\,k_v\,|\,w\in W,\,v\in V\}$ is called {\it the labelled
graph} of $M$ and is denoted by $X_M$ (here and further we omit
``arguments" of $X_M$).

 \subheading{1.4 The BKN-equation} Let $M\in\M$ be a
 graph manifold and let $X_M$
 be its labelled graph. The equation
 $$
 \sum_{w\in\, \partial v}\frac{\gamma_wa_{w(v)}}{b_{w}}=k_v\cdot a_v,
 $$
 (the symbol $w(v)$ denotes the terminal
 vertex for $w$)
  with $\{\,a,\,\gamma\,|\,v\in V,\,w\in W\}$
 as unknown rational numbers
 is called the {\it BKN-equation} over $X_M$.
The \BKN is said to be {\it solvable}
 if there exists a nontrivial ($a\not\equiv 0$)
 rational solution $\{a,\gamma\}$ such that
 $a_v\ge 0$, $|\gamma_{w}|\le 1$ and $\gamma_w\gamma_{-w}\ne -1$.

 The following lemma describes the correspondence between CP-systems
 and solutions of the BKN-equation.

 \proclaim{1.5 Lemma} Let $M$ be a graph manifold of the class $\M$.
 The \BKN over $X_M$
 is solvable if and only if
 there exists a \CP $\{\,l_v\,|\,{v\in V}\}$ on $M$.
 \endproclaim

 \demo{Proof}
 Let $\{a,\,\gamma\}$ be a solution of the \BKN, and
 let
$$
W_0=\{\,w\in W\,|\,\text{if}\;w\in\partial v \;\text{then}\;
 a_va_{w(v)}=0\}.
$$
Now we define a new solution $\{a,\,\gamma'\}$ of the \BKN by
$$
\gamma'_w=\left\{\matrix
\gamma_w&w\not\in W_0\\
0&w\in W_0
\endmatrix\right..
$$

Consider an element
 l$_w\in H^1(T_{|w|})$ such that
 $$
 \text{l}_w(f_w)=a_v,\quad
 \text{l}_w(f_{-w})=\gamma_w'a_{w(v)},
 $$
 and the element l$_v\in H^1(\partial M_v)$
 such that $p_w($l$_v)=$l$_w$ for each $w\in\partial v$.
 Using lemma 1.3.1 it is not difficult to verify that the element
 l$_v$ lies in the image of $H^1(M_v)$ under the canonical map $\iota^*$
  which is induced by inclusion $\iota:\partial M_v\to M_v$.
 Now we get the required
 $l_v$ as an (arbitrary) element of $(\iota^*)^{-1}\text{l}_v$.
 If l$_v=0$ we choose $l_v=0$.

 Conversely, let $\{\,l_v\,|\,{v\in V}\}$ be a \CP on $M$.
 To find a solution of the \BKN we can use the above two formulas:
 $$
 a_v=l_v(f_v)\quad\text{ and}\quad
 \gamma_w=\left\{\matrix
 \frac{\langle l_v,f_{v'}\rangle_{|w|}}{a_{w(v)}}& \text{if}
 \,\,\, a_{w(v)}\ne 0, \\
 &\\
 0 & \text{otherwise}.
 \endmatrix
 \right.
 $$

 So we get
 $$
 \multline
 \sum_{w\in\, \partial v}\frac{\gamma_wa_{w(v)}}{b_{w}}=
 \sum_{w\in\, \partial v}\frac{\langle l_v,f_{w(v)}\rangle_{|w|}}{b_{w}}=
 \sum_{w\in\, \partial v}
 \iota^*_wl_v\left(\frac{f_{-w}}{b_{w}}\right) \\
 =l_v\left(\sum_{w\in\, \partial v}(\iota_w)_*\frac{f_{-w}}{b_{w}}\right)
 =-e_{M_v}(\phi_v)l_v(f_v)=k_va_v.
 \qed
 \endmultline
 $$
 \enddemo

 It is easy to translate theorem I from CP-language
 to BKN-language.
 \proclaim{Theorem II}~A manifold $M\in\M$ satisfies the
 property {\bf Im} (resp. {\bf HI}, {\bf F}, {\bf  E}, {\bf VF}, {\bf  VE},
 {\bf NPC})
 if and only if the \BKN over $\Gamma_M(V,W)$
 has a solution (resp. has a solution $\{a,\,\gamma\}$ such that
  \roster
 \item "(HI)" $a_v>0$, and if $|\gamma_w|=1$ then $\gamma_{-w}=\gamma_{w}$;
 \item "(F)" $a_v>0$ and $\gamma_w=\gamma_{-w}=\epsilon_{v'}\epsilon_v$
 for a function $\epsilon:V\to\{\pm 1\}$;
 \item "(E)" either 1) $\gamma_w=\gamma_{-w}=\pm 1$ or
 2) $a_v=0$ and $\gamma_w=\gamma_{-w}=0$ for any $w\in\partial v$;
 \item "(VF)" $a_v>0$  and $\gamma_w=\gamma_{-w}$;
 \item "(VE)" $\gamma_w=\gamma_{-w}$;
 \item "(NPC)" $a_v>0$ and $\gamma_w=\gamma_{-w}\in (-1;1)$,
 \endroster
 for each $v\in V$, $w\in W$).
 \endproclaim
 This theorem follows from theorem I and lemma 1.5. \qed
 \subheading{Remark} Items {\bf Im} and {\bf HI} of theorem II
 was proved by W.~Neumann in \cite{N2}.
 The item {\bf NPC} was proved by S.~Buyalo and V.~Kobel'skii in \cite{BK1}.

\subheading{1.6 Matrices and reduced criteria} Consider the
following square matrices $A^{\epsilon}_M=(a^{\epsilon}_{vv'})$,
$A^+_M=(a^+_{vv'})$, and $H_M=(h_{vv'})$ over the graph
$\Gamma_M$:
$$
a^{\epsilon}_{vv'}=\left\{\matrix
k_v-\sum\limits_{w(v)=v}\frac{\epsilon_w}{b_w},& \text{if}~ v=v'\\
&\\
-\sum\limits_{w(v)=v'}\frac{\epsilon_w}{b_w},&
 \text{if}~v\ne v'
\endmatrix\right.,
\quad
a^+_{vv'}=\left\{\matrix
|k_v|-\sum\limits_{w(v)=v}\frac{1}{|b_w|},& \text{if}~ v=v'\\
&\\
-\sum\limits_{w(v)=v'}\frac{1}{|b_w|},&
 \text{if}~v\ne v'
\endmatrix\right.,
$$

$$
h_{vv'}=\left\{\matrix
s(v)k_v-\sum\limits_{w(v)=v}\frac{1}{|b_w|},&\text{if}~ v=v'\phantom{GGGGGGG}\\
&\\
-\sum\limits_{w(v)=v'}\frac{1}{|b_w|},
& \text{if}~ v\ne v'\text{and}~ k_v\cdot k_{v'}>0\\
&\\
0,&\text{otherwise}\phantom{GGGGGG}
\endmatrix\right..
$$
Here $\epsilon\:W\to\{\pm 1\}$ is a function.
The matrix $-A_M^{\epsilon}$ for $\epsilon\equiv 1$ is the
``reduced plumbing matrix" from \cite{N1}.
The second matrix is equal to $-A_-$ from \cite{N2}.
The function $s:V\to \{0,\pm 1\}$ is constructed as follows.
Vertices $v,v'\in V$
of the graph $\Gamma_M$ are called equivalent $v\sim v'$
if there exists a path $v_0=v$, $v_1$, $v_2$, $\dots$, $v_n=v'$
in the graph $\Gamma_M$ such
that $k_{v_i}\cdot k_{v_{i+1}}>0$, for each $i=0,\dots,\,n-1$.
An edge $w\in W$ and a vertex $v\in V$ are called equivalent $v\sim w$
if $k_v\ne 0$ and $v$ is equivalent to both the ends of $w$.
Edges $w,w'\in W$ are called equivalent if either  $w'=-w$ or there exists
a vertex $v\in V$ such that $w\sim v\sim w'$.
Now we define the graph of signed components $G(U,E_0)$
of the triple $\{\Gamma_M(V,W), B, K\}$ as the factor graph $\Gamma_M/\sim$.
This graph $G(U,E_0)$ was firstly defined in \cite{BK2}.

If the graph $G(U, E_0)$ is not bipartite one then we put
$s(v)=0$ for each vertex  $v\in V$. Let $G(U,E_0)$ be
a bipartite graph and $U=P\cup N$ be a splitting to parts such that
$P\cap\{\,v\in V\,|\,k_v>0\}\ne\emptyset$. In this case we put
$s(v)=1$ if $v\in P$ and $s(v)=-1$ if $v\in N$.

If the graph $G(U,E_0)$ is bipartite then the matrix $H_M$
 coincides with a matrix from \cite{BK2}.
 If a graph manifold $M$
 has no block with zero charge then $G(U;E_0)$ is bipartite
 and $s(v)=\sgn k_v$, if in addition the graph $\Gamma_M$
has no loops then  $H_M$ can be represented as
 the matrix $-(P_-\oplus N)$ from the paper \cite{N1}.

A matrix $A$ is called {\it semipositive defined} (resp. {\it supersingular})
if $x^tAx\ge 0$ (resp. $Ax=0$) for
each tuple $x$ (resp. for some tuple $x$ with no zero entry).
A square matrix $A'$ is called a {\it principal submatrix} of a square
matrix $A$ if it can be obtained from $A$ by deleting
some rows and the corresponding columns.

\proclaim{Theorem III}
A manifold $M\in\M$ satisfies the
property {\bf Im}, {\bf HI}, {\bf F}, {\bf  E}, {\bf VF}, {\bf  VE},
and {\bf NPC}
if and only if the following conditions (respectively) hold
 \roster
\item "(Im, HI)" either 1) diagonal elements of $A_M$
has the same sign and the matrix $A_M^+$ is semipositive defined and singular or
2) the matrix $A_M^+$ has a negative eigenvalue;
\item "(F)" the matrix $A_M$ is supersingular;
\item "(E)" there exists $\epsilon\: W\to\{\pm 1\}$ such that
the matrix $A^{\epsilon}_M$ has a  supersingular principal submatrix;
\item "(VF)" either 1) the matrix $H_M$ is semipositive defined
 and  supersingular
or 2) the matrix $H_M$ has a negative eigenvalue;
\item "(VE)"the matrix $H_M$ has a  principal submatrix
which satisfies the requirement of the previous item
\item "(NPC)" either 1) $H_M\equiv 0$
or 2) the matrix $H_M$ has a negative eigenvalue.
\endroster
\endproclaim

The item {\bf F} of this theorem was proved in \cite{N1},
the items {\bf Im} and {\bf VE}  were proved in \cite{N1,2}
for manifolds whose graphs has no loops,
the item {\bf NPC} was proved in \cite{BK2}, and the
item {\bf VF} was proved by the author in \cite{S}.
We give no proof of the remain item {\bf E} in this draft version.\qed

\subheading{1.7 Conventions}
Let $g:S\to M$ be a $\pi_1$-injective immersion.
Rubinstein and Wang \cite{RW} proved that each
surface $S$ with negative Euler characteristic in a graph-manifold $M$
is  $\pi_1$-injective if and only
if it can be properly homotoped so that any connected component
of $S\cap M_v$
is either vertical or horizontal for each maximal block $M_v$.
Further we assume that all immersed surfaces already put in this
position.

Let $g:S\to M$ be a $\pi_1$-injective not horizontal immersion
and let $g':S'\to M$ be
the immersion such that $g'(S')$ is the boundary of a collar
for $g(S)$ in $M$.
\comment
Let $A_i\subset g'(S')\setminus{\Cal T},\,i=1,2$ be parallel annuli
in a maximal block $M_v$ which join tori
$T_1\subset \partial M_v$ and $T_2\subset \partial M_v$.
Let $A_i'$ be not boundary parallel annuli
in $M_v$ such that
$\partial A'_i=\partial(A_1\cup A_2)\cap T_i$ for each $i=1,2$.
Let $g'':S''\to M$ be an immersion such that
$g''(S'')=\left(g'(S')\setminus (A_1\cup A_2)\right)\cup A'_1\cup A_2'.$
\endcomment
If $g:S\to M$ is a $\pi_1-$injective immersion
(embedding) then $g':S'\to M$
has the same property.
If $g:S\to M$ is a $\pi_1-$injective
virtual embedding  then $g':S'\to M$
has the same property. Converse is also true \cite{RW}.

Further we assume that each vertical annulus of $g(S)\setminus{\Cal T}$
has a parallel copy.

 \head 2.  Immersed and horizontal immersed surfaces in graph manifolds
\endhead

\subheading {2.1}  The implication i$\Leftrightarrow$iii
of the following theorem was proved by W.~Neumann in \cite{N2}.
\proclaim{Theorem (cf. [N3])}Let $M$ be a
manifold of the class $\M$. The following three properties
are equivalent.

\roster
\item "i"  $M$ contains an immersed $\pi_1$-injective surface
of negative Euler characteristic.

\item "ii" there exists a \CP  on $M$

\item "iii" \BKN over $\Gamma_M$ is solvable.
   \endroster
\endproclaim

\demo{Proof}~ Let $g:S\to M$ be a $\pi_1$-injective immersion. and
let ${\Cal C}=g^{-1}({\Cal T})\subset S$. The set $S\setminus
{\Cal C}$ is a disjoint union $\cup_{\alpha\in A}S_{\alpha}$ of
connected components. Now we orient the surfaces
$S_{\alpha},\,\alpha\in A$ in the following way. Choosing
orientation of $M$ and orientations of the Seifert fibers in
maximal blocks induce an orientation on each surface $S_{\alpha}$
which is not annulus. By conventions 1.7, vertical part of
$g(S)\cap M_v$ is a set of parallel pairs of annuli. We orient the
annuli so that parallel annuli have opposite orientations.

Let $T_{|w|}\subset\partial M_v\cap\partial M_{v'}$ be a
JSJ-torus and $C$ be a connected component of $g^{-1}(T_{|w|})$.
The curve $C$ is oriented as a boundary component
of $g^{-1}(M_v)$  and at the same time it is oriented
as a boundary component
of $g^{-1}(M_{v'})$.
The curve is said to be {\it consistent}
(resp. {\it nonconsistent}) if 
these two orientations
coincide (resp. opposite).
The sum of homological classes of consistent
(resp. nonconsistent) curves on $T_{|w|}$ which are oriented
as boundary components of $g(S)\cap M_v$
we denote by $c_w^+$ (resp. $c_w^-$), where $w\in\partial v$.

By Poincare duality $H_2(M_v;\partial M_v)\simeq H^1(M_v)$
we get a class $l_v\in H^1(M_v)$ which is dual for
$[g(S)\cap M_v]\in H_2(M_v;\partial M_v)$.
It is clear that
$$
\iota_w^*l_v(x)=x\wedge_w(c_w^++c_w^-),\quad
\iota_{-w}^*l_{v'}(x)=x\wedge_w(c_w^+-c_w^-)
$$
for each $x\in H_1(T_{|w|})$, where
$w\in\partial v$, $-w\in\partial v'$.
Due to orientations on components of $g(S)\cap M_v$
we have $a_w^+=f_w\wedge_wc_w^+\ge 0$ and $a_w^-=f_w\wedge_wc_w^-\ge 0
$.

Since
$$
l_v(f_v) = a_w^++a_w^-
,\quad
\langle l_{v'},f_v\rangle_{|w|}= a_w^+-a_w^-
$$
we have
$|\langle l_{v'},f_v\rangle_{|w|}|\le l_v(f_v)$ (CP1).

If either $a_w^+=0$ or $a_w^-=0$
then (by convention 1.7) we conclude that no one of curves
$g(S)\cap T_{|w|}$ is parallel to Seifert fiber of $M_v$
as well as of $M_{v'}$.
Therefore
either $a_w^-=a_{-w}^-=0$ or $a_w^+=a_{-w}^+=0$.
So we have
$\langle l_v,\cdot\rangle_{|w|}=\pm \langle l_{v'},\cdot\rangle_{|w|}$
and it is equal to CP2.

Finally, if $l_v(f_v)=0$ then $g(S)\cap M_v$ is a set of pairs
of oppositely oriented parallel annuli. So we have $l_v=0$ (CP3).

ii$\Leftrightarrow$iii It is the assertion of lemma 1.5.

ii$\Rightarrow$i~Multipluing all $\{l_v\}_{v\in V}$ by
suitable integer, we can assume that $l_v$
is an integer class for each $v\in V$.

Let $c_w^{+}\in H_1(T_{|w|})$
(resp. $c_w^-\in H_1(T_{|w|})$) be the homological class
which is dual by Poincare for
$
\iota_w^*l_v+\iota_{-w}^*l_{v'}\in H^1(T_{|w|})$
(resp. for  $\iota_w^*l_v-\iota_{-w}^*l_{v'}\in H^1(T_{|w|})$), where
$w\in\partial v,\,-w\in\partial v'$.
It is obvious that $c_{-w}^{\pm}=\pm c_w^+$.
Choose  a pair of
nonoriented curves
$C_{|w|}^+$ and $C_{|w|}^-$ on each JSJ-torus $T_{|w|}$ in $M$
so that
these curves represent (with some orientation)
the classes $c_w^+$ and $c_w^-$ respectively.
If one of the classes (or both) is zero then
we do not consider the corresponding curve.

 It is easy to see that
$$
f_w\wedge_w(c_w^{+}+c_w^{-})=2l_v\left((\iota_w)_*f_w\right)=2l_v(f_v),
\tag{$1'$}$$
$$
\sum_{w\in\partial v}p_w(z_v)\wedge_w(c_w^{+}+c_w^{-})
=2l_v\left(\sum_{w\in\partial v}(\iota_w)_*p_w(z_v)\right)=
-2e_{M_v}(z_v)\cdot l_v(f_v).
\tag{$2'$}
$$
for each adjoint element $z_v\in H_1(\partial M_v),\,v\in V$.
If $p_w(z_v)=f_{-w}/b_w$ for each $w\in \partial v$
then the right side of \thetag{$2'$} is $k_vl_v(f_v)$.
The equalities \thetag{$1'$} and \thetag{$2'$} coincide
with the corresponded equalities from lemma 1.3.1 if $l_v(f_v)>0$.

As above, by $d(C)^{n}$ we denote $d$ copies of a curve obtained
by going $n$ times around the curve $c$.
Using lemma 1.3.1 we can find
integers $d_v,\,n_v$ such that for any integer pair $d,\, n$
with $d_v$ dividing $d$ and  $n_v$ dividing $n$
the curves
$d(C^{\pm}_{|w|})^{n},\,w\in\partial v$
span a surface (horizontal if $l_v(f_v)\ne 0$)
in $M_v$. Taking appropriate integers $D$ and $N$ we
obtain a surface by fitting together the parts $S_v$ spanning on
$D(C_{|w|}^{\pm})^N,\,w\in\partial v$ for each $v\in V$
such that $l_v(f_v)>0$.
If $l_v(f_v)>0$ for each $v\in V$ we obtain the desired surface.

Consider a block $M_v\subset M$ such that $l_v=0$.
For each $T_w\subset\partial M_v$ we either have
two set of parallel curves $D(C_{|w|}^+)^N$ and $D(C_{|w|}^-)^N$
on $T_{|w|}$ or have no curves.
In the first case the curves are parallel to
the linear foliation of $T_w$ which is induced from $M_v$.
So we can find D parallel immersed vertical annuli in $M_v$ which are not
boundary parallel and their boundary is
$D(C_{|w|}^+)^N\cup D(C_{|w|}^-)^N$.
The horizontal parts $S_v$ (in blocks $M_v$ with $l_v(f_v)\ne 0$) as above
and the annuli fit together providing
the desired surface.\qed
\enddemo

\subheading{2.2~ The HI property}
In this section we prove
\proclaim{Theorem (cf. [N2])}Let $M$ be a
manifold of the class $\M$. The following three properties
are equivalent.

\roster
\item "i" $M$ contains a horizontal immersed surface
of negative Euler characteristic.

\item "ii" There is a \CP $\{\,l_v\,|\,v\in V\}$ on M such that
$l_v(f_v)>0$ for each $v\in V$.

\item "iii" \BKN over $\Gamma_M$ has a nontrivial solution $\{a,\gamma\}$
such that $a_v>0$ for each $v\in V$.
   \endroster
\endproclaim

\demo{Proof}
The equivalence ii$\Leftrightarrow $iii follows from lemma 1.5.

$[i\Rightarrow ii]$~Using previous theorem it sufficient to prove
the condition $l_v(f_v)>0$.
It is true, due to horizontality.

$[ii\Rightarrow i]$ ~The surface fitting as in the
last part of the previous proof is horizontal in
each block $M_v$ inasmuch $l_v(f_v)>0$
for each $v\in V $ and the surface is connected in each maximal block.
\qed
\enddemo


 \head          3.  Fibering and embedded surfaces   \endhead

\comment
~Let $M$ be a three-manifold which is fibered over $S^1$
with surface $S$ as a fiber.
We can obtain the manifold identifying bases of the cylinder
$S\times [0;\,1]$ by a homeomorphism $\phi:S\to S$ as follows:
$(x,0)\sim(\phi (x),1)$.

The strong Thurston result \cite{Th, Theorem 0.1} implies
that if the above $M$
is a graph manifold of the class $\M$ then
\roster
\item $\phi $ is isotopic to a reducible homeomorphism $\ov{\phi}$
\item the invariant set ${\Cal C}$ of $\ov{\phi}$ is
${\Cal T}\cap (S\times \{pt\})$ --  the union of traces
of the JSJ-tori on a fiber
\item the restriction of $\phi$ to each connected component
of $S\setminus {\Cal C}$ is isotopic to an
irreducible homeomorphism of finite degree.
\endroster
\endcomment
The map $p:M\to S^1$ is homotopic to a fibration iff
$p_*(f_v)\ne 0$ for each $v\in V$ \cite{NR}.
So, if $p$ is a fibration then any fiber $p^{-1}(pt)$
(surface)
is transversal to Seifert fibers (circles) in each block
so such a surface is {\it horizontal}.

\proclaim{3.1 Theorem (cf. \cite{N1})}
Let $M$ be a manifold of the class $\M$ and $\Gamma_M(V,W)$ be its graph.
The following conditions are equivalent.
\roster
\item "i" $M$ is fibered over $S^1$ with oriented surface
of negative Euler characteristic as a fiber.

\item "ii" There exists a \CP $\{\,l_v\,\}_{v\in V}$  on $M$
such that $l_{v}(f_{v})>0$
for each $v\in V$ and in addition for each triple $v,w,v'$
as in CP1 (definition 1.2) we have
$\langle l_{v'},f_{v}\rangle_{|w|}=
\epsilon_{v'}\epsilon_{v}l_{v}(f_{v})$, where
$\epsilon:V\to\{\pm 1\}$ is a function.

\item "iii" There exists a solution $\{a,\gamma\}$ of \BKN
such that $a_v>0$ for each $v\in V$
and $\gamma_w=\gamma_{-w}=\epsilon_{v'}\epsilon_{v}$ for each triple
$v,w,v'$ as in CP1 (definition 1.2), where
$\epsilon:V\to\{\pm 1\}$ is a function.

\endroster

            \endproclaim

\demo{Proof}
The equivalence ii$\Leftrightarrow $iii follows from lemma 1.5.

i$\Rightarrow$ii
Let $p\:M\to S^1$ be a fibering map and $p^{-1}(pt)$
is a fiber $S$, $\chi(S)<0$.
For each $v\in V$ by $p_v:M_v\to S^1$ we denote the restriction
of $p$ to the maximal block $M_v$.

Put by definition $\epsilon_v=\sgn\langle \alpha,p_*f_v\rangle$
where $\alpha $ is a generator of $H_1(S^1;\Bbb Z)$.
The desired classes $l_v\in H^1(M_v),\,v\in V$ we define as follows
$$
l_v=\epsilon_v\cdot p_v^*\alpha,
$$
It is clear that $l_v(f_v) =|\langle \alpha,p_*f_v\rangle|>0$
and is equal to the geometric intersection
number between $S\cap M_v$ and a Seifert fiber in $M_v$
which are transversal each other. To conclude the proof
of $i\Rightarrow ii$ it remains to note that

$$
\langle l_v,f_{v'}\rangle_{|w|}=
\epsilon_v\langle p_v^*\alpha,f_{v'}\rangle_{|w|}=
\epsilon_v\langle \iota_w^*p_v^*\alpha,f_{-w}\rangle=
\epsilon_v\langle \alpha,p_*f_{v'}\rangle=
\epsilon_v\epsilon_{v'}l_{v'}(f_{v'}).
$$

[$ii\Rightarrow i$]
Due to $\epsilon_vl_v(f_v)\ne 0$
we can find a fibering map
$p_v:M_v\to S^1$ using the canonical isomorphism
$H^1(M_v)\simeq [M_v,S^1]$ so that $p_v^*(\alpha)=\epsilon_vl_v$ \cite{N1}.
Since $\iota_w^*\epsilon_vl_v=\iota_{-w}^*\epsilon_{v'}l_{v'}$
for each triple $v,w,v'$
we can choose fibering maps $p_v$ and $p_{v'}$
so that $p_v|_{T_{|w|}}=p_{v'}|_{T_{|w|}}$
for each $T_{|w|}\subset M_v\cap M_{v'}$.
By \cite{NR} there is a fibering map $p\:M\to S^1$
such that $p|_{M_v}=p_v$.
\qed
\enddemo

\proclaim{3.2 Theorem (cf. \cite{N1})}
Let $M$ be a manifold of the class $\M$ and $\Gamma_M(V,W)$ be its graph.
The following conditions are equivalent.
\roster
\item "i" $M$ contains an embedded surface
of negative Euler characteristic.

\item "ii" There exists a \CP $\{l_v\}_{v\in V}$
on  $M$ such that
for each triple $v,w,v'$
as in CP1 (definition 1.2) we have
$|\langle l_{v'},f_{v}\rangle_{|w|}|=l_v(f_v)$.

\item "iii" There exists a solution $\{a,\gamma\}$ of the \BKN
over $\Gamma_M(V,W)$
such that either $\gamma_w=\gamma_{-w}=\pm 1$
or $a_v=0$ and $\gamma_w=\gamma_{-w}=0$ for each $w\in\partial v$.

\endroster
\endproclaim

\demo{Proof}~
The equivalence ii$\Leftrightarrow$iii follows from
lemma 1.5.

[i$\Rightarrow$ii] Let $g\:S\to M$ be a $\pi_1$-injective
embedding. Consider the boundary  $S'$ of $N(g(S))$, where $N$
denotes a normal neighborhood.
 We orient the connected components of $S'\setminus{\Cal T}$ as
in the proof of theorem 2.1. Let $l_v\in H^1(M_v)$ be the class, which is dual
for $[S'\cap M_v]\in H_2(M_v,\partial M_v)$. Since
the surface $S'$ has negative Euler characteristic, we have
$l_{v^*}(f_{v^*})\ne 0$ at least for one block $M_{v^*}\subset M$.
Consider a JSJ-torus $T_{|w|}$ separating blocks $M_v$ and
$M_{v'}$. Since $g$ is an embedding, each intersection $S'\cap M_v$
either vertical or horizontal. If $T_{|w|}$ does not intersect the
embedded surface then $\iota_w^*l_v=\pm \iota_{-w}^*l_{v'}=0$
and we have nothing to do.
Assume that $l_v(f_v)>0$
then either $l_{v'}(f_{v'})>0$ or
$l_{v'}(f_{v'})=0$.
In the first case we have $\iota_w^*l_v=\pm \iota_{-w}^*l_{v'}$
indeed $M_v$ and $M_{v'}$
contain only horizontal parts of the {\it embedded} surface
so we get
$\langle l_{v'},f_{v}\rangle_{|w|}=\pm l_{v}(f_{v})$.
In the second case pairs of opposite oriented annuli are occurred in
$M_{v'}$
and corresponded pairs of parallel horizontal
surfaces are occurred in $M_v$. It gives
$\langle l_{v'},f_{v}\rangle_{|w|}=\langle l_{v},f_{v'}\rangle_{|w|}=0$
therefore $\iota_{-w}^*l_{v'}=0$ for each $-w\in\partial v'$.
In any case we have $|\langle l_{v'},f_{v}\rangle_e|=l_v(f_v)$
and $|\langle l_{v},f_{v'}\rangle_e|=l_{v'}(f_{v'})$.
So $\{\,l_v\,|\,v\in V\}$ is the desired \CP.

[ii$\Rightarrow$i]~
Let $\{\,l_v\,|\,v\in V\}$ be a \CP which is satisfied ii
and $T_{|w|}$ be a JSJ-torus in $M$ such that $w\in\partial v$, and
$-w\in\partial v'$.
As in the theorem's 2.1 proof we can consider
the pair of homological classes $c_w^+$ and $c_w^-$.
There are three cases:
\comment
so that the sum of its homological classes $c_w^+,\,c_w^-\in H_1(T_{|w|})$
be dual to $\iota^*_wl_v$ on $T_{|w|}$.                       \endcomment
\roster
\item both $l_v$ and $l_{v'}$ are not zero,
\item either $l_v=0$ or $l_{v'}=0$ but not both,
\item both $l_v$ and $l_{v'}$ are zero.
\endroster
In the first case ii and CP2 imply $\iota_{-w}^*l_{v'}=\pm\iota_{w}^*l_v$
and there is only one non-zero homological class (either $c_w^+$ or $c_w^-$).
In the second case either $\iota_{-w}^*l_{v'}(f_w)=0$ or
$\iota_{w}^*l_{v}(f_{-w})=0$ respectively
and we have two parallel curves which are parallel to the linear
foliation on $T_{|w|}$ which is induced either
from $M_v$ or from $M_{v'}$.
In the third case there are no curves on $T_{|w|}$.

By lemma 1.3.2 there exists an embedded $\pi_1$-injective surface
$S\subset M$ such that the curves in $S\cap{\Cal T}$ realize the
classes $c_w^+$, $c_w^-$ for all $w\in W$. \qed
\enddemo

 \head          4. Virtual fibering and
virtually embedded surfaces
\endhead
Let $M$ be a graph manifold of the class $\M$. A covering
$p:\widetilde{M}\to M$ is called {\it s-cha\-rac\-te\-ris\-tic}
(or characteristic), if its restriction on each JSJ-torus
${T}\subset\widetilde{M}$ is organized as follows: there exists a
basis $(\tilde{a},\tilde{b})$ of the group $\pi_1({T})$ and a
basis $({a},{b})$ of the group $\pi_1(p\,({T}))$ such that
$(p|_{{T}})_*\tilde{a}=sa$ and $(p|_{{T}})_*\tilde{b}=sb$ for a
positive integer $s$ \cite{LW}. Let
$q:\Gamma_{\widetilde{M}}(\widetilde{V},\widetilde{W})\to
\Gamma_M(V,W)$ -- be the map between the graphs which is induced
by $p$, and let $d_us^2$  be the degree of the restriction
$p|_{M_u}:M_u\to M_{q(u)},\,u\in\widetilde{V}.$ Note, that the
number $s^2\sum_{q(u)=v}d_u$ is independent of $v\in V$ and is
equal to the degree of $p$. In the following lemma we collect some
useful results from \cite{LW}. 

\proclaim{4.1 Lemma (Luecke and Wu)} 
Let $p:\widetilde{M}\to M$ be a finite covering of graph manifold $M$. 
\roster
\item There exists a covering $\tilde{p}:\widehat{M}\to\widetilde{M}$
such that  $\hat{p}=p\circ\tilde{p}$
is a characteristic covering \cite{LW, prop.4.4}. In this case
\item the charges of vertices and the intersection indexes for edges
of $\Gamma_{\widehat{M}}$ и $\Gamma_M$
are connect by the following relations:
$k_u=k_{\hat{q}(u)}d_u$, $b_l=b_{\hat{q}(l)}$,
where $u\in\widehat{V},\,l\in\widehat{W}$ and
$\hat{q}:\Gamma_{\widehat{M}}\to\Gamma_M$ is the map induced
by $\hat{p}$ \cite{LW, prop.2.3}.
\endroster
\endproclaim

 \proclaim{4.2 Proposition}
 Let $M$ be
a graph manifold of the class $\M$ and let $\Gamma_M(V,W)$ be
 its graph. If $M$ contains a virtually embedded
$\pi_1$-injective surface of negative Euler
characteristic then \BKN has a solution $\{a,\gamma\}$
such that  $\gamma_w=\gamma_{-w}$ for each $w\in W.$
Moreover, if the surface is horizontal
then $a_v>0$ for each $v\in V$.
 \endproclaim

\demo{Proof}
Let $M$ satisfies {\bf VE}.
By item 1 of the previous lemma we have
an $s$-cha\-rac\-te\-ri\-stic covering
$p:\widehat{M}\to M$ and a $\pi_1-$injective embedding
$g:S\to \widehat{M}.$
By theorem 3.2 there exists a solution $\{x,\mu\}$
of the \BKN over $\Gamma_{\widehat{M}}$
such that $\mu_l=\mu_{-l}\in\{0,\pm 1\}$ for each $l\in\widehat{W}$:
$$
k_ux_u=\sum_{l\in\partial u}\frac{\mu_l}{b_l}\,x_{l^+},
$$
where $l^+$ is the terminal edge for $l$.
Using the item 2 of lemma 4.1 we have
$$
k_{\hat{q}(u)}d_ux_u=\sum_{w\in\partial \hat{q}(u)}\frac{1}{b_w}
\sum_{l\in\partial u\cap\hat{q}^{-1}(w)}{\mu_l}x_{l^+}.
\tag{*}
$$

Let $M_w,\;w\in {W}$ be a rectangular matrix given by
$$
\left(M_w\right)_{uu'}=
\sum_{\{l\in \hat{q}^{-1}(w)|\,u'=l(u)\}}\mu_l
$$
The columns of $M_w$ are corresponded to vertices of
$\hat{q}^{-1}(w^+)$ as well as its rows are corresponded to
vertices of $\hat{q}^{-1}(w^-)$. Since $\mu_l=\mu_{-l}$ we have
$M_{-w}=M_w^t$, where $t$ is the transposition. Let $X_v,\;v\in V$
be a tuple consisting of entries $x_u,\;u\in \hat{q}^{-1}(v)$ and
let $D_v$ be a diagonal matrix $\diag\{d_u,\;u\in
\hat{q}^{-1}(v)\}.$ In these notations the equality \thetag{*} can
be rewrited as a vector one:
$$
k_v D_vX_v=\sum_{w\in \partial v} \frac{1}{b_w}M_wX_{w^+}.
\tag{**}
$$

Put $t_v=+\sqrt{X_v^tD_vX_v}$ and
$$
\lambda_w=
\frac{X_{w^-}^tM_wX_{w^+}}
{t_{w^-}\cdot t_{w^+}}\quad \text{if}\quad t_{w^-}\cdot t_{w^+}\ne 0
$$
and $\lambda_w=0,$ {if} $t_{w^-}\cdot t_{w^+}=0.$
Multiplying both sides of \thetag{**} by $X_{w^-}^t$
from the left, we get
  $$
  k_vt_v^2=\sum_{w\in\partial v}\frac{\lambda_w}{b_{w}}\,t_vt_{w^+}.
  $$
Since $\lambda_w=0$ for each $w\in\partial v$
{if} $t_{v}=0,$
we have
  $$
  k_vt_v=\sum_{w\in\partial v}\frac{\lambda_w}{b_{w}}\,t_{w^+}.
  $$
\enddemo

\subheading{Claim} $|\lambda_w|\le 1$.

If $A=(a_{ik})$ is $m\times n$-matrix,
$x=(x_i)$ -- $m$-tuple,
$y=(y_k)$ -- $n$-tuple, and all entries
$a_{ik},\;x_i,\;y_k$ are nonnegative then
$$
\multline
x^tAy=\sum_{i,k}x_ia_{ik}y_k=
\sum_ix_i\left(\sum_k\sqrt{a_{ik}}\cdot\sqrt{a_{ik}}y_k\right)\le\\
\sum_i\left(x_i\sqrt{\sum_ka_{ik}}\right)\cdot
\sqrt{\sum_ka_{ik}y_k^2}\le\\
\sqrt{\sum_i\left(\sum_ka_{ik}\right)x_i^2}\cdot
\sqrt{\sum_k\left(\sum_ia_{ik}\right)y^2_k}
\endmultline
$$
by  the Cauchy inequality acting twice.

Let $M'_w$ be the matrix which is obtained from $M_w$
by taking absolute values of all its entries.
It is clear that $|X^t_vM_wX_{w(v)}|\le X^t_vM'_wX_{w(v)}.$
The previous inequality with the matrix $M'_w$
and the tuples $X_v,\;X_{w(v)},$ gives
$$
X^t_vM'_wX_{w(v)}\le
\sqrt{X^t_vD^{|w|}_vX_v\cdot X^t_{w(v)}D^{|w|}_{w(v)}X_{w(v)}},
$$
where $D^{|w|}_v$ is the diagonal matrix with diagonal entries
which is defined by
$$d^{|w|}_u=\sum_{\sssize l\in\partial u\cap \hat{q}^{-1}(w)}
|\mu_l|,\quad u\in\hat{q}^{-1}(v).$$
Since $|\mu_l|\le 1$ we have
$d^{|w|}_u\le \#\{\partial u\cap \hat{q}^{-1}(w)\}=d_u$
hence
$$
\multline
|X^t_vM_wX_{w(v)}|\le
X^t_vM'_wX_{w(v)}\le
\sqrt{X^t_vD^{|w|}_vX_v\cdot X^t_{w(v)}D^{|w|}_{w(v)}X_{w(v)}}\le\\
\sqrt{X_v^tD_vX_v\cdot X_{w(v)}^tD_{w(v)}X_{w(v)}}=t_vt_{w(v)}.
\endmultline
$$
 The claim is proved.

Note, that  $\lambda_{-w}=\lambda_w$ since
$M_w^t=M_{-w}$.
If the surface $g(S)$ is horizontal then
$t_v>0$ for each $v\in V$.
So we have a solution $\{t,\lambda\}$ of the \BKN  over
$\Gamma(V,W),$
with prescribed  properties.
\qed

\proclaim{4.3 Theorem}Let $M$ be a graph manifold
of the class $\M$. The following three properties
are equivalent.

\roster
\item "i"  $M$ is virtually fibered over a circle.

\item "ii" There exists a \CP of cohomological classes
$\{\,l_v\,|\,v\in V\}$ on $M$
such that $l_v(f_v)> 0$ for each $v\in V$, and
$\langle l_{v'},f_{v}\rangle_{|w|}\cdot l_{v'}(f_{v'})=
\langle l_{v},f_{v'}\rangle_{|w|}\cdot l_{v}(f_{v})$
for each triple $v,w,v'$ as in CP1 (definition 1.2).

\item "iii" The \BKN over $\Gamma_M$
has a nontrivial solution $\{a,\gamma\}$
such that $a_v>0$ and $\gamma_w=\gamma_{-w}$ for each $v\in V$, $w\in W$.
   \endroster
\endproclaim

\demo{Proof}
The equivalence ii$\Leftrightarrow$iii follows from lemma 1.5.
The implication i$\Rightarrow$ii is in proposition 4.2.

The implication  ii$\Rightarrow$i~ is proved in our previous
work \cite{S}. The proof uses the main result of \cite{RW}.\qed
\enddemo

\comment
I have a little time to take myself away.

I have no evident in my unknown mission.

As you could see, I spend all gifted blessom, wisdom.

My cup is emptied not by me.

\endcomment

Now we need some technique which is useful in what follows.
Let $M$ be a manifold of the class $\M$, let $\Gamma_M(V,W)$ be its graph
and let $V_0\subset V$ be a subset of its maximal blocks.
Remove the blocks $\{M_v|\,v\in V_0\}$ from $M$ and let
$$
M\setminus V_0:=\overline{M\setminus\bigcup_{v\in V_0}M_v}
$$
be a (possible not connected) manifold with boundary. The manifold
$M\setminus V_0$ is canonically framed as follows. Consider
$T_{|w|}$, a boundary torus of $M\setminus V_0$ (here the initial
point of $w$ is in $V_0$ and its terminal point is in $V\setminus
V_0$), and choose $r_w$ as a curve on $T_w$ which corresponds to
the Seifert fibration of the removing block $M_{v}$, $w\in\partial
v$. Let $(M\setminus V_0)^D$ be the closed manifolds by performing
Dehn filling along all $r_w$, $T_{|w|}\subset\partial (M\setminus
V_0)$. This manifold is defined up to homeomorphism. The charges
of maximal blocks of $({M\setminus V_0})^D$ are coincides with the
charges of corresponding blocks of $M$ \cite{N1}. The graph of
$({M\setminus V_0})^D$ is a subgraph of $\Gamma_M(V,W).$

\proclaim{4.4 Theorem (cf. \cite{N1})}
Let $M$ be a manifold of the class $\M$ and $\Gamma_M(V,W)$ be its graph.
The following conditions are equivalent.
\roster
\item "i"  $M$ contains a virtually embedded $\pi_1$-injective
surface of negative Euler characteristic.
\item "ii" There exists a \CP of cohomological classes
$\{\,l_v\,|\,v\in V\}$ on $M$
such that
$\langle l_{v'},f_{v}\rangle_{|w|}\cdot l_{v'}(f_{v'})=
\langle l_{v},f_{v'}\rangle_{|w|}\cdot l_{v}(f_{v})$
for each triple $v,w,v'$ as in CP1 (definition 1.2).

\item "iii" The \BKN over $\Gamma_M$
has a solution $\{a,\gamma\}$
such that $\gamma_w=\gamma_{-w}$ for each $v\in V$, $w\in W$.
   \endroster
\endproclaim

\demo{Proof}~ The equivalence ii$\Leftrightarrow$iii follows from lemma 1.5.

The implication i$\Rightarrow$ii is in proposition 4.2.

iii$\Rightarrow$i~ Let $\{a,\gamma\}$ be a solution of the
 \BKN which satisfies iii and $V_0=\{v\in V\,|\, a_v=0\}$.
Consider a manifold $M_1$ with boundary which is a connected part
of ${M\setminus V_0}$, and let $M_1^D$ be the corresponding closed
graph manifold (may be Seifert fibered). It is easy to see that
the collection $\{a_v,\gamma_w|v\in V_1, w\in W_1\}$ is a solution
of the \BKN over $\Gamma_{M_1}(V_1,W_1)\subset \Gamma_{M}(V,W)$
and it satisfies to iii of theorem 4.3. By the theorem, there
exists a finite covering $p^D_1:\widetilde{M_1^D}\to M_1^D$ and a
horizontal embedding
$\widetilde{g^D_1}:S^D_1\to\widetilde{M^D_1}$, where $S^D_1$ is a
closed oriented surface.
 Let $p_1:\widetilde{M_1}\to M_1$ be the covering which
is induced by the inclusion $\iota_1:M_1\to M_1^D$, i.e.
$p_1=\iota_1^*p_1^D$. We also have a natural proper horizontal
embedding $\widetilde{g_1}\:S_1\to \widetilde{M_1}$, where $S_1$
is the surface $S_1^D$ with some number of holes and
$\chi(S_1)<0$.

Now we claim that there exists a covering $p:\widetilde{M}\to M$
such that $\iota^*(p)=p_1$ for the natural inclusion $\iota\:M_1\to M$.
Really, this inclusion induces the injective map
$\iota_*:\pi_1(M_1)\to \pi_1(M)$.
 Such a way we have an embedding
$(S_1,\partial S_1)\to(\widetilde{M},\widetilde{{\Cal T})}$
where $\widetilde{{\Cal T}}$ is the JSJ-surface for $\widetilde{M}$.
The surface $S_1$ is contained in $M_1\subset M$.
Let $C$ be a connected component of $\partial S_1\cap T$
where $T\in \widetilde{{\Cal T}}$ and $\partial S\cap T\ne\emptyset$.
The torus $T$ must separates two blocks. One of the blocks say $M_u$
must be preimage (under $p$) of one of removing blocks $M_v$, $v\in V_0$.
As it easy to see the curve $C$ is homotopic to the linear
foliation of $T$ which is induced from the Seifert space $M_u$.
Let $(S_2,\partial S_2)\subset (\widetilde{M},\widetilde{{\Cal T})}$
be a surface obtained as two parallel copies of $S_1$.
The intersection $M_u\cap S_2$ consists on even number of curves on
$\partial M_u$. Each curve parallel to the linear
foliation of the boundary component it lies.
It is not difficult to see that there are a set of incompressible,
boundary-incomp\-ressible annuli in $M_u$ whose boundary is $M_u\cap S_2$.
So we have found a $\pi_1$-injective surface in some finite cover of $M$.
 \qed \enddemo

 \head          5. NPC-metrics on graph manifolds
\endhead

\proclaim{5.1 Theorem (cf. \cite{BK2})}
Let $M$ be a manifold of the class $\M$ and $\Gamma_M(V,W)$ be its graph.
The following conditions are equivalent.
\roster
\item "i" $M$
carries an NPC-metric.
\item "ii" there exists a \CP $\{l_v\}_{v\in V}$ on $M$ such that
|$\langle l_{v'},f_{v}\rangle_e|< l_{v}(f_{v})$
and
$\langle l_{v'},f_{v}\rangle_e\cdot l_{v'}(f_{v'})=
\langle l_{v},f_{v'}\rangle_{|w|}
\cdot l_{v}(f_{v})$
for each triple $v,e,v'$
as in CP1 (definition 1.2);

\item "iii" there exists a solution $\{a,\gamma\}$ of CE
such that $a_v>0$ for each $v\in V$
and $\gamma_w=\gamma_{-w}\in (-1,1)$ for each $w\in W$.

\endroster
\endproclaim

\demo{Proof}~
The equivalence ii$\Leftrightarrow$iii follows from lemma 1.5.

 \,iii$\Leftrightarrow$i~
If each Seifert block of $M$ is a product surface$\times S^1$
then the equivalence follows from \cite{BK1, Proposition 8.1}.
The general case is proved by I.~Andreeva (unpublished).
 \qed
 \enddemo

\comment

 \head          6. one criterion
\endhead

\proclaim{6.1 Proposition}
A graph manifold $M\in\M$
satisfies to the property {\bf E}
if and only if
there exists $\epsilon\: W\to\{\pm 1\}$ such that
the matrix $A^{\epsilon}_M$ has a  supersingular principal submatrix.
\endproclaim

\demo{Proof}
Let $g:S\to M$ be a $\pi_1$-injective embedding
and let $V_0$ be a set of blocks with vertical annuli.
Consider a connected component $M_1$
of $(M\setminus V_0)^D$ (this notation is explained before theorem 4.4).
Let $\Gamma_{M_1}(V_1,W_1)\subset\Gamma_M(V,W)$ be the corresponded subgraph
and let $\{a,\,\gamma\}$ be the
solution of the \BKN over $\Gamma_M(V,W)$ (by theorem 4.4).
В силу того, что $a_v=0$ для любой вершины $v\in V\setminus{V_1},$
соседствующей с какой-либо вершиной из $V_1,$
набор $\{a_v,\,\gamma_w\,|\,v\in V_1\,w\in W_1\}$
является решением БКН-уравнения над $\Gamma_{M_1}(V_1,W_1),$
при этом $a_v>0$ и $|\gamma_w|=1.$
Возьмем многообразие $M_1'$ 2-соизмеримое (п. 1.6.3) c $M_1$ такое,
что $b_w'=\gamma_wb_w$ для всех $w\in W_1.$
Набор $\{a_v,\,\gamma_w\,|\,v\in V_1\,w\in W_1\}$
доставляет решение БКН-уравнения над раскрашенным графом многообразия
$M_1',$ удовлетворяющий требованиям п. {\bf F} теоремы II.
Согласно этой теореме $M_1'$ расслоено над окружностью,
а значит матрица $A_{M_1'}$ (соизмеримая с $A_M$),
(супер)сингулярна по предложению 3.2.
Осталось заметить, что $A_{M_1}$ является
главной подматрицей матрицы $A_M.$

Обратно, если у $A_M$ имеется
главная подматрица, соизмеримая с сингулярной,
то, у нее есть и некоторая неприводимая
 главная подматрица $A_1,$
соизмеримая с суперсингулярной.
Для этого достаточно показать, что у сингулярной матрицы $A$ есть
суперсингулярная главная подматрица $B.$
Пусть $x, y, \cdots,z$ --- линейно независимые стобцы,
порождающие ядро матрицы $A.$ Если матрица $A$
не суперсингулярна, то $x_i=y_i=\cdots=z_i=0$ для
$i-$ых компонент этих столбцов при некотором $i.$
Пусть $A_1$ --- матрица, полученная вычеркиванием $i-$го столбца и
соответствующей строки из матрицы
$A.$ Эта матрица опять особая.
Применим описанную выше процедуру к матрице $A_1,$ и т.д.
За конечное число шагов получим искомую матрицу $B.$
(Для матриц размером $1\times 1$ сингулярность и суперсингулярность
эквивалентны.)

Пусть $\Gamma(V_1,W_1)$ --- (связный) подграф графа многообразия $M,$
над которым определена матрица $A_1,$ соизмеримая
с некоторой суперсингулярной матрицей $A_1',$
а $(l_v)_{v\in V_1}$ --- столбец из ядра матрицы $A_1'$
без нулевых элементов. Положим
$$
a_v=\left\{\matrix
|l_v|,&v\in V_1\\
0,&v\not\in V_1
\endmatrix\right.,\quad
\gamma_w=\left\{\matrix
\sgn(l_{w^-}l_{w^+}b_wb'_{w}),&w\in W_1\\
0,&w\not\in W_1
\endmatrix\right..
$$
Этот набор удовлетворяет условию iii теоремы 2.4 для многообразия $M$,
поэтому оно обладает свойством {\bf E}.
\qed
\enddemo
\endcomment

$$$$

e-mail {svetlov\@pdmi.ras.ru} 

\end{document}